\documentclass{article}

\title{General self-similarity: an overview}
\author{Tom Leinster}
\date{\normalsize 
Department of Mathematics, University of Glasgow\\
www\dt maths\dt gla\dt ac\dt uk\slsh $\sim$tl\\
tl@$\,\!$maths\dt gla\dt ac\dt uk}

\input{diagrams}
\diagramstyle[repositionpullbacks,midshaft,PostScript=dvips,nohug,scriptlabels]
\usepackage{graphics}
\usepackage[mathcal,mathscr]{euscript}

\usepackage{amssymb}

\newcommand{\dashbk}{-}
\newcommand{\diagspace}{\mbox{\hspace{2em}}}
\newcommand{\bref}[1]{(\ref{#1})}

\newcommand{\mb}[1]{\mathbf{#1}}
\newcommand{\mc}[1]{\mathcal{#1}}
\newcommand{\mr}[1]{\mathrm{#1}}
\newcommand{\cat}[1]{\mc{#1}}
\newcommand{\fcat}[1]{\mb{#1}}

\newcommand{\slsh}{/\linebreak[0]}
\newcommand{\dt}{.\linebreak[0]}

\newcommand{\epsln}{\varepsilon}

\newcommand{\op}{\mr{op}}

\newcommand{\Set}{\fcat{Set}}
\newcommand{\Top}{\fcat{Top}}
\newcommand{\ftrcat}[2]{[#1,#2]}
\newcommand{\goesto}{\,\longmapsto\,}

\newcommand{\parpair}[2]{\pile{\rTo^{\scriptstyle #1}\\ 
\rTo_{\scriptstyle #2}}}
\newcommand{\parpairu}{\pile{\rTo\\ \rTo}}
\newarrow{Get}....>
\newarrow{Goesto}|---{->}
\newarrow{Incl}C--->
\newarrow{Mod}--+->

\newcommand{\latin}[1]{\textit{#1}}
\newcommand{\demph}[1]{\textbf{\textup{#1}}}
\newcommand{\done}{\hfill\ensuremath{\Box}}
\newenvironment{prooflike}[1]{\begin{trivlist}\item\textbf{#1}\ }
{\end{trivlist}}
\newenvironment{proof}{\begin{prooflike}{Proof}}{\end{prooflike}}
\newcommand{\scat}[1]{\mathbb{#1}}
\newcommand{\go}{\rTo\linebreak[0]}

\newcommand{\iso}{\cong}
\newcommand{\nat}{\mathbb{N}}	
\newcommand{\url}[1]{#1}

\newcommand{\pshf}[1]{\ftrcat{#1^\op}{\Set}}

\newcommand{\ndSet}[1]{[#1, \Set]_{\mr{nondegen}}}

\newcommand{\cstyle}[1]{\textbf{\upshape{#1}}}
\newcommand{\So}{\cstyle{S}}

\newcommand{\dfrc}[2]{\textstyle{\frac{#1}{#2}}}
\newcommand{\dhalf}{\dfrc{1}{2}}
\newcounter{bean}
\newcommand{\Cx}{\mathbb{C}}

\newcommand{\cell}[4]{\put(#1,#2){\makebox(0,0)[#3]{\ensuremath{#4}}}}
\newcommand{\numcell}[3]{\put(#1,#2){\makebox(0,0){\ensuremath{\scriptstyle
#3}}}}

\newtheorem{thm}{Theorem} 

\newtheorem{lemma}[thm]{Lemma}
\newtheorem{conj}{Conjecture}
\newtheorem{defn}{Definition}

\begin{document}
\sloppy

\maketitle

\begin{abstract}
Consider a self-similar space $X$.  A typical situation is that $X$ looks
like several copies of itself glued to several copies of another space $Y$,
and that $Y$ looks like several copies of itself glued to several copies of
$X$---or the same kind of thing with more than two spaces.  Thus, we have a
system of simultaneous equations in which the right-hand sides (the
gluing instructions) are `higher-dimensional formulas'.

This idea is developed in detail in~\cite{SS1} and~\cite{SS2}.  The present
informal seminar notes explain the theory in outline.
\end{abstract}

\noindent
I want to tell you about a very general theory of self-similar objects that
I've been developing recently.  In principle this theory can handle
self-similar objects of any kind whatsoever---algebraic, analytic,
geometric, probabilistic, and so on.  At present it's the topological case
that I understand best, so that's what I'll concentrate on today.  This
concerns the self-similar or fractal spaces that you've all seen pictures
of.

Some of the most important self-similar spaces in mathematics are Julia
sets.  For the purposes of this talk you don't need to know the definition
of Julia set, but the bare facts are these: to every complex rational
function $f$ there is associated a closed subset $J(f)$ of the Riemann
sphere $\Cx \cup \{\infty\}$, its Julia set, which almost always has a
highly intricate, fractal appearance.  If you look in a textbook on complex
dynamics, you'll find theorems about `local self-similarity' of Julia sets.
For example, given almost any point $z$ in a Julia set, points locally
isomorphic to $z$ occur densely throughout the set \cite[Ch.~4]{Mil}.  On
the other hand, the kind of self-similarity I'm going to talk about today
is the dual idea, `global self-similarity', where you say that the
\emph{whole} space looks like several copies of itself stuck together---or
some statement of the kind.  So it's a top-down, rather than bottom-up,
point of view.

A long-term goal is to develop the algebraic topology of self-similar spaces.
The usual invariants coming from homotopy and homology are pretty much
useless (e.g.\ for a fractal subset of the plane all you've got is $\pi_1$,
which is usually either trivial or infinite-dimensional), but a description
by global self-similarity is discrete and so might provide useful
invariants at some point in the future.

This theory is about self-similarity as an \emph{intrinsic} structure on an
object: there is no reference to an ambient space, and in fact no ambient
space at all.  This is like doing group theory rather than representation
theory, or the theory of abstract manifolds rather than the theory of
manifolds embedded in $\mathbb{R}^n$.  We can worry about representations
later.  For instance, the Koch snowflake is just a circle for us: its
self-similar aspect is the way it's embedded in the plane.

Later I'll show you the general language of self-similarity, but first here
are some concrete examples to indicate the kind of situation that I want to
describe.

\section{First example: a Julia set}

\begin{figure}
\setlength{\unitlength}{1mm}
\begin{picture}(120,50)
\cell{21}{50.5}{t}{\includegraphics{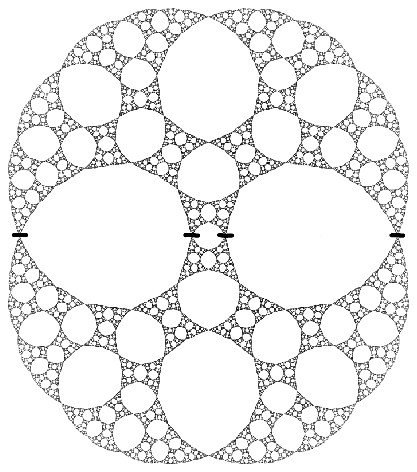}}
\cell{65}{7}{b}{\includegraphics{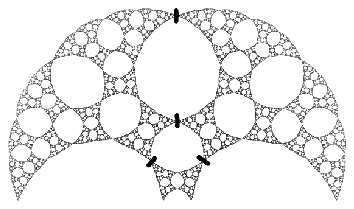}}
\cell{104}{7}{b}{\includegraphics{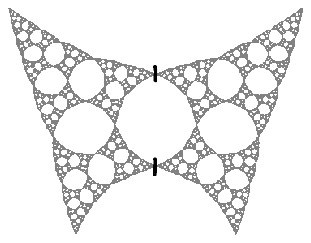}}
\cell{21}{-2}{b}{\textrm{(a)}}
\cell{65}{-2}{b}{\textrm{(b)}}
\cell{49}{6.5}{t}{1}
\cell{63.2}{6.5}{t}{2}
\cell{66.8}{6.5}{t}{3}
\cell{81}{6.5}{t}{4}
\cell{104}{-2}{b}{\textrm{(c)}}
\cell{96}{6.5}{t}{1}
\cell{88.5}{31}{b}{2}
\cell{119}{31}{b}{3}
\cell{112}{6.5}{t}{4}
\end{picture}
\caption{(a) The Julia set of $z \goesto (2z/(1 + z^2))^2$; (b), (c)
  certain subsets}
\label{fig:Julia}
\end{figure}

Let's look at one particular Julia set in detail: Figure~\ref{fig:Julia}(a).  
I'll write $I_1$ for this Julia set.  It clearly has reflectional symmetry
in a horizontal axis, so if we cut at the four points shown then we get a
decomposition
\begin{equation}	\label{eq:I1}
I_1 = 
\begin{array}{c}
\setlength{\unitlength}{1mm} 
\begin{picture}(18,20)(-9,-10)
\cell{0.1}{0.1}{c}{\includegraphics{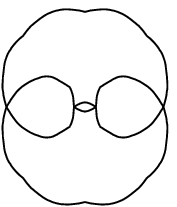}}
\numcell{-7.7}{2}{1}
\numcell{-0.9}{1.8}{2}
\numcell{0.9}{1.8}{3}
\numcell{7.6}{2}{4}
\numcell{-7.7}{-2}{1}
\numcell{-0.9}{-1.8}{2}
\numcell{0.9}{-1.8}{3}
\numcell{7.6}{-2}{4}
\cell{0}{6}{c}{I_2}
\cell{0}{-6}{c}{I_2}
\end{picture}
\end{array}
\end{equation}
where $I_2$ is a certain space with four marked points (or `basepoints').
Now consider $I_2$ (Figure~\ref{fig:Julia}(b)).  Cutting at the points
shown gives a decomposition
\begin{equation}	\label{eq:I2}
I_2 
=
\begin{array}{c}
\setlength{\unitlength}{1mm} 
\begin{picture}(52,30)(-26,0)
\cell{0.2}{0}{b}{\includegraphics{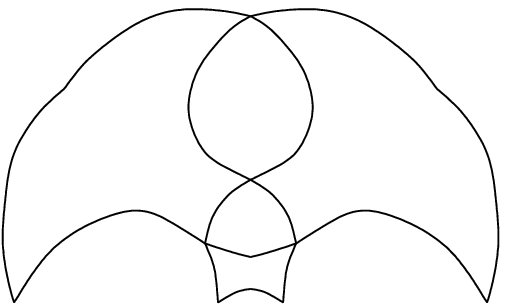}}
\numcell{-23.6}{2.5}{1}
\numcell{-5}{7.8}{2}
\numcell{-2.3}{12.5}{3}
\numcell{-3}{28.5}{4}
\numcell{23.8}{2.5}{1}
\numcell{5.2}{7.8}{2}
\numcell{2.5}{12.5}{3}
\numcell{3}{28.5}{4}
\numcell{-2.8}{1.7}{1}
\numcell{-3.3}{4.8}{2}
\numcell{3.3}{4.8}{3}
\numcell{2.8}{1.7}{4}
\cell{-13}{16}{c}{I_2}
\cell{13}{16}{c}{I_2}
\cell{0}{2}{b}{I_3}
\end{picture}
\end{array}
\end{equation}
where $I_3$ is another space with four marked points.  Next, consider $I_3$
(Figure~\ref{fig:Julia}(c)); it decomposes as
\begin{equation}	\label{eq:I3}
I_3 
= 
\begin{array}{c}
\setlength{\unitlength}{1mm} 
\begin{picture}(32,20)(-16,0)
\cell{0}{0}{b}{\includegraphics{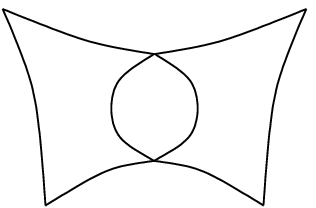}}
\numcell{-3}{5.5}{1}
\numcell{-10}{2.5}{2}
\numcell{-13}{18}{3}
\numcell{-3}{15}{4}
\numcell{3}{5.5}{1}
\numcell{10}{2.5}{2}
\numcell{13}{18}{3}
\numcell{3}{15}{4}
\cell{-8}{10}{c}{I_3}
\cell{8}{10}{c}{I_3}
\end{picture}
\end{array}.
\end{equation}
Here we can stop, since no new spaces are involved.  Or nearly: there's a
hidden role being played by the one-point space $I_0$, since that's what we've
been gluing along, and I'll record the trivial equation
\begin{equation}	\label{eq:I0}
I_0 = I_0.
\end{equation}

What we have here is a system of four simultaneous equations, with the
unusual feature that the right-hand sides are not algebraic formulas of the
usual type, but rather `2-dimensional formulas' expressing how the spaces
are glued together.  

(There's a conceptual link between this and the world of $n$-categories,
where there are 2-dimensional and higher-dimensional morphisms which you're
allowed to compose or `glue' in various ways.  Both can be regarded as a
kind of higher-dimensional algebra.  The \latin{cognoscenti} will see a
technological link too: in both contexts the gluing can be described by
pullback-preserving functors on categories of presheaves.)

The simultaneous equations \bref{eq:I1}--\bref{eq:I0} can be expressed as
follows.  First, we have our spaces $I_1, I_2, I_3$ with their marked
points, which together form a functor from the category
\[
\scat{A} = 
\left( \ \ 
\begin{diagram}[height=1.7em,width=3em,noPS]
	&	&1	\\
	&	&	\\
0	&\pile{\rTo\\ \rTo\\ \rTo\\ \rTo}	
		&2	\\
	&\pile{\rdTo\\ \rdTo\\ \rdTo\\ \rdTo}
		&	\\
	&	&3	\\
\end{diagram}
\ \ 
\right)
\]
to the category $\Set$ of sets (or a category of spaces, but
let's be conservative for the moment).  Second, the gluing
formulas define a functor
\[
G: \ftrcat{\scat{A}}{\Set} \go \ftrcat{\scat{A}}{\Set},
\]
where $\ftrcat{\scat{A}}{\Set}$ is the category of functors $\scat{A} \go
\Set$: given $X \in \ftrcat{\scat{A}}{\Set}$, put 
\begin{eqnarray*}
(G(X))_1	&=	&
\begin{array}{c}
\setlength{\unitlength}{1mm} 
\begin{picture}(18,20)(-9,-10)
\cell{0.1}{0.1}{c}{\includegraphics{X1new.eps}}
\numcell{-7.7}{2}{1}
\numcell{-0.9}{1.8}{2}
\numcell{0.9}{1.8}{3}
\numcell{7.6}{2}{4}
\numcell{-7.7}{-2}{1}
\numcell{-0.9}{-1.8}{2}
\numcell{0.9}{-1.8}{3}
\numcell{7.6}{-2}{4}
\cell{0}{6}{c}{X_2}
\cell{0}{-6}{c}{X_2}
\end{picture}
\end{array}
=
(X_2 \amalg X_2)/\sim,	
\\
(G(X))_2	&=	&
\begin{array}{c}
\setlength{\unitlength}{1mm} 
\begin{picture}(52,31)(-26,0)
\cell{0.2}{0}{b}{\includegraphics{X2new.eps}}
\numcell{-23.6}{2.5}{1}
\numcell{-5}{7.8}{2}
\numcell{-2.3}{12.5}{3}
\numcell{-3}{28.5}{4}
\numcell{23.8}{2.5}{1}
\numcell{5.2}{7.8}{2}
\numcell{2.5}{12.5}{3}
\numcell{3}{28.5}{4}
\numcell{-2.8}{1.7}{1}
\numcell{-3.3}{4.8}{2}
\numcell{3.3}{4.8}{3}
\numcell{2.8}{1.7}{4}
\cell{-13}{16}{c}{X_2}
\cell{13}{16}{c}{X_2}
\cell{0}{2}{b}{X_3}
\end{picture}
\end{array}
=
(X_2 \amalg X_2 \amalg X_3)/\sim ,
\end{eqnarray*}
and so on.  (I've drawn the pictures as if $X_0$ were a single point.)
Then the simultaneous equations assert precisely that
\[
I \iso G(I)
\]
---$I$ is a fixed point of $G$.

Although these simultaneous equations have many solutions ($G$ has many
fixed points), $I$ is in some sense the universal one.  This means that the
simple diagrams \bref{eq:I1}--\bref{eq:I0} contain just as much information
as the apparently very complex spaces in Figure~\ref{fig:Julia}: for given
the system of equations, we recover these spaces as the universal solution.
Caveats: we're only interested in the intrinsic, topological aspects of
self-similar spaces, not how they're embedded into an ambient space (in
this case, the plane) or the metrics on them.

Next we have to find some general way of making rigorous the idea of
`gluing formula'; so far I've just drawn pictures.  We have a small
category $\scat{A}$ whose objects index the spaces involved, and I claim
that the self-similarity equations are described by a functor $M:
\scat{A}^\op \times \scat{A} \go \Set$ (a `2-sided $\scat{A}$-module').
The idea is that for $b, a \in \scat{A}$,
\begin{eqnarray*}
M(b, a)		&=	& 
\{\textrm{copies of the }b\textrm{th space used in the gluing formula} \\
	&	&
\ \,\textrm{for the }a\textrm{th space} \}.
\end{eqnarray*}
Take, for instance, our Julia set.  In the gluing formula for $I_2$, the
one-point space $I_0$ appears 8 times, $I_1$ doesn't appear at all, $I_2$
appears twice, and $I_3$ appears once, so, writing $n$ for an
$n$-element set,
\[
M(0, 2) = 8,
\diagspace
M(1, 2) = \emptyset,
\diagspace
M(2, 2) = 2,
\diagspace
M(3, 2) = 1.
\]
(It's easy to get confused between the arrows $b \go a$ in $\scat{A}$ and
the elements of $M(b, a)$.  The arrows of $\scat{A}$ say nothing whatsoever
about the gluing formulas, although they determine where gluing may
\emph{potentially} take place.  The elements of $M$ embody the gluing
formulas themselves.)

This is an example of what I'll call a `self-similarity system':
\begin{defn} \label{defn:sss}
A \demph{self-similarity system} is a small category $\scat{A}$ together
with a functor $M: \scat{A}^\op \times \scat{A} \go \Set$ such that
\begin{enumerate}
\item \label{item:Mfinite}
for each $a \in \scat{A}$, the set $\coprod_{c, b \in \scat{A}} \scat{A}(c,
b) \times M(b, a)$ is finite
\item \label{item:Mpbf}
(a condition to be described later).
\end{enumerate}
\end{defn}
Part~\bref{item:Mfinite} says that in the system of simultaneous equations,
each right-hand side is a gluing of only a \emph{finite} family of spaces.
So we might have infinitely many spaces (in which case $\scat{A}$ would be
infinite), but each one is described as a finite gluing.  The condition is
more gracefully expressed in categorical language: `for each $a$, the
category of elements of $M(\dashbk, a)$ is finite'.

As in our example, any self-similarity system $(\scat{A}, M)$ induces an
endofunctor $G$ of $\ftrcat{\scat{A}}{\Set}$.  This works as follows.
First note that if $A$ is a ring (not necessarily commutative), $Y$ a right
$A$-module, and $X$ a left $A$-module, there is a tensor product $Y
\otimes_A X$ (a mere abelian group).  Similarly, if $\scat{A}$ is a small
category, $Y: \scat{A}^\op \go \Set$ a contravariant functor, and $X:
\scat{A} \go \Set$ a covariant functor, there is a tensor product
\[
Y \otimes X 
=
\left( \coprod_{a \in \scat{A}} Ya \times Xa \right) / \sim
\]
(a mere set): see~\cite[IX.6]{CWM}.  So if $(\scat{A}, M)$ is a
self-similarity system then there is an endofunctor $G = M \otimes \dashbk$
of $\ftrcat{\scat{A}}{\Set}$ defined by
\[
(M \otimes X)(a) 
= 
M(\dashbk, a) \otimes X
=
\left( \coprod_{b \in \scat{A}} M(b, a) \times Xb \right) / \sim.
\]
($X \in \ftrcat{\scat{A}}{\Set}$, $a \in \scat{A}$).  We are interested in
finding a fixed point of $G$ that is in some sense `universal'.

\section{Second example: Freyd's Theorem}

The second example I'll show you comes from a very different direction.  In
December 1999, Peter Freyd posted a message~\cite{Fre} on the categories
mailing list that caused a lot of excitement, especially among the
theoretical computer scientists.

We'll need some terminology.  Given a category $\cat{C}$ and an endofunctor
$G$ of $\cat{C}$, a \demph{$G$-coalgebra} is an object $X$ of $\cat{C}$
together with a map $\xi: X \go GX$.  (For instance, if $\cat{C}$ is a
category of modules and $GX = X \otimes X$ then a $G$-coalgebra is a
coalgebra---not necessarily coassociative---in the usual sense.)  A
\demph{map} $(X, \xi) \go (X', \xi')$ of coalgebras is a map $X \go X'$ in
$\cat{C}$ making the evident square commute.  Depending on what $G$ is, the
category of $G$-coalgebras may or may not have a terminal object, but if it
does then it's a fixed point:
\begin{lemma}[Lambek~\cite{Lam}]
Let $\cat{C}$ be a category and $G$ an endofunctor of $\cat{C}$.  If $(I,
\iota)$ is terminal in the category of $G$-coalgebras then $\iota: I \go
GI$ is an isomorphism.
\end{lemma}
\begin{proof}
Short and elementary.
\done
\end{proof}

Here's what Freyd said, modified slightly.  Let $\cat{C}$ be the category
whose objects are diagrams $X_0 \parpair{u}{v} X_1$ where $X_0$ and $X_1$
are sets and $u$ and $v$ are injections with disjoint images; then an
object of $\cat{C}$ can be drawn as
\[
\setlength{\unitlength}{1mm}
\begin{picture}(30,14)(-15,-5)
\cell{0}{0}{c}{\includegraphics{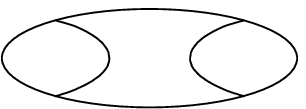}}
\cell{-10}{0}{c}{X_0}
\cell{10}{0}{c}{X_0}
\cell{0}{5.5}{b}{X_1}
\end{picture}
\]
where the copies of $X_0$ on the left and the right are the
images of $u$ and $v$ respectively.  A map $X \go X'$ in
$\cat{C}$ consists of functions $X_0 \go X'_0$ and $X_1 \go X'_1$
making the evident two squares commute.  Now, given $X \in
\cat{C}$ we can form a new object $GX$ of $\cat{C}$ by gluing
two copies of $X$ end to end:
\[
\begin{array}{c}
\setlength{\unitlength}{1mm}
\begin{picture}(49,14)(-24.5,-5)
\cell{0}{0}{c}{\includegraphics{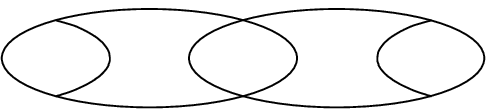}}
\cell{0}{0}{c}{X_0}
\cell{-20}{0}{c}{X_0}
\cell{20}{0}{c}{X_0}
\cell{-10}{5.5}{b}{X_1}
\cell{10}{5.5}{b}{X_1}
\end{picture}
\end{array}.
\]
Formally, $GX$ is defined by pushout:
\[
\begin{diagram}[size=2em]
	&	&	&	&(GX)_1	&	&	&	&	\\
	&	&	&\ruTo	&	&\luTo	&	&	&	\\
	&	&X_1	&	&\textrm{pushout} &
						&X_1	&	&	\\
	&\ruTo<u&	&\luTo>v&	&\ruTo<u&	&\luTo>v&	\\
(GX)_0 = X_0&	&	&	&X_0	&	&	&	&X_0.	\\
\end{diagram}
\]
For example, the unit interval with its endpoints distinguished forms an
object 
\[
I = \left( \{\star\} \parpair{0}{1} [0,1] \right)
\]
of $\cat{C}$, and $GI$ is naturally described as an interval of
length 2, again with its endpoints distinguished:
\[
GI = \left( \{\star\} \parpair{0}{2} [0,2] \right).
\]
So there is a coalgebra structure $\iota: I \go GI$ on $I$ given by
multiplication by two.  Freyd's Theorem says that this is, in fact, the
\emph{universal} example of a $G$-coalgebra:
\begin{thm}[Freyd$+\epsln$]	\label{thm:Freyd}
$(I, \iota)$ is terminal in the category of $G$-coalgebras.
\end{thm}

I won't go into the proof, but clearly it's going to have to involve the
completeness of the real numbers in an essential way.  Once you've worked
out what `terminal coalgebra' is saying, it's easy to see that the proof is
going to be something to do with binary expansions.  Note that although
$\iota$ is an isomorphism (as predicted by Lambek's Lemma), this by no
means determines $(I, \iota)$: consider, for instance, the unique coalgebra
satisfying $X_0 = X_1 = \emptyset$, or the evident coalgebra in which $X_0 =
\{\star\}$ and $X_1 = [0,1] \cap \{ \textrm{dyadic rationals} \}$.

The striking thing about Freyd's result is that we started with just some
extremely primitive notions of set, function, and gluing---and suddenly,
out popped the real numbers.  What excited the computer scientists was that
it suggested new ways of representing the reals.  But its relevance for us
here is that it describes the self-similarity of the unit interval---in
other words, the fact that it's isomorphic to two copies of itself stuck
end to end.  We'll take this idea of Freyd, describing a very simple
self-similar space as a terminal coalgebra, and generalize it dramatically.

Freyd's Theorem concerns the self-similarity system $(\scat{A}, M)$ in
which
\[
\scat{A} = 
\left( 0 \parpairu 1 \right)
\]
and $M: \scat{A}^\op \times \scat{A} \go \Set$ is given by
\[
\begin{diagram}
M(\dashbk, 0):&	&\{\star\}	&\pile{\lTo\\ \lTo}	&\emptyset	\\
	&	&\dTo<0 \dTo>1	&			&\dTo \dTo	\\
M(\dashbk, 1):&	&\{ 0, \dhalf, 1\}&
				\pile{\lTo^{\inf}\\ \lTo_{\sup}}	&
\{ [0, \dhalf], [\dhalf, 1] \}.	\\
\end{diagram}
\]
(Here $M(0, 1)$ is just a $3$-element set and $M(1, 1)$ a $2$-element set,
but their elements have been named suggestively.)  The category $\cat{C}$
is a full subcategory of $\ftrcat{\scat{A}}{\Set}$, and the endofunctor $G$
of $\cat{C}$ is the restriction of the endofunctor $M \otimes \dashbk$ of
$\ftrcat{\scat{A}}{\Set}$.

The only thing that looks like a barrier to generalization is the condition
that $u, v: X_0 \go X_1$ are injective with disjoint images (which is the
difference between $\cat{C}$ and $\ftrcat{\scat{A}}{\Set}$).  If this were
dropped then $(\{\star\} \parpairu \{\star\})$ would be the terminal
coalgebra, so the theorem would degenerate entirely.  It turns out that the
condition is really a kind of flatness.

A module $X$ over a ring is called `flat' if the functor $\dashbk \otimes
X$ preserves finite limits.  There is an analogous definition when $X$ is a
functor, but actually we want something weaker:
\begin{defn}
Let $\scat{A}$ be a small category.  A functor $X: \scat{A} \go \Set$ is
\demph{nondegenerate} if the functor
\[
\dashbk \otimes X: \pshf{\scat{A}} \go \Set
\]
preserves finite \emph{connected} limits.  Write $\ndSet{\scat{A}}$ for the
category of nondegenerate functors $\scat{A} \go \Set$ and natural
transformations between them.
\end{defn}
It looks as if this is very abstract, that in order to show that $X$ was
degenerate you'd need to test it against all possible finite connected
limits in $\pshf{\scat{A}}$, but in fact there's an equivalent explicit
condition.  This can be used to show that in the case at hand, a functor
$X: \scat{A} \go \Set$ is degenerate precisely when the two functions $u,
v: X_0 \go X_1$ are injective with disjoint images.

The missing condition~\bref{item:Mpbf} in the definition of self-similarity
system is that for each $b \in \scat{A}$, the functor $M(b, \dashbk):
\scat{A} \go \Set$ is nondegenerate.  This guarantees that the endofunctor
$M \otimes \dashbk$ of $\ftrcat{\scat{A}}{\Set}$ restricts to an
endofunctor of $\ndSet{\scat{A}}$.  A terminal coalgebra for this
restricted endofunctor is called a \demph{universal solution} of the
self-similarity system; if one exists, it's unique up to canonical
isomorphism.  Lambek's Lemma implies that if $(I, \iota)$ is a universal
solution then, as the terminology suggests, $M \otimes I \iso I$.  Freyd's
Theorem describes the universal solution of a certain self-similarity
system.

Before we move on, I'll show you a version of Freyd's Theorem in which the
unit interval is characterized not just as a set but as a topological
space.  The simplest thing would be to change `set' to `space' and
`function' to `continuous map' throughout the above.  Unsurprisingly, this
gives a boring topology on $[0,1]$: the indiscrete one, as it happens.  But
all we need to do to get the usual topology is to insist that the maps $u$
and $v$ are closed.  That is:
\setcounter{bean}{\value{thm}}
\begin{trivlist}
\item
\textbf{Theorem \arabic{bean}$'$}
\textit{Define $\cat{C'}$ and $G'$ as $\cat{C}$ and
  $G$ were defined above, changing `set' to `space' and `function' to
  `continuous map', and adding the condition that $u$ and $v$ are closed
  maps.  Then the terminal $G'$-coalgebra is $(I,\iota)$ where $[0,1]$ is
  equipped with the Euclidean topology.}
\end{trivlist}
Generally, a functor $X: \scat{A} \go \Top$ is \demph{nondegenerate} if its
underlying $\Set$-valued functor is nondegenerate and $Xf$ is a closed map
for every map $f$ in $\scat{A}$.  This gives a notion of \demph{universal
topological solution}, just as in the set-theoretic scenario.  So
Theorem~\arabic{bean}$'$ describes the universal topological solution of
the Freyd self-similarity system.

\section{Results}

Just as some systems of equations have no solution, some self-similarity
systems have no universal solution.  But it's easy to tell
whether there is one:
\begin{thm} 
  A self-similarity system has a universal solution if and only if it
  satisfies a certain condition \So.
\end{thm}
I won't say what \So\ is, but it is completely explicit.  So too is the
construction of the universal solution when it does exist; it is similar
in spirit to constructing the real numbers as infinite decimals, although
smoother than that would suggest.

Let $(\scat{A}, M)$ be a self-similarity system with universal solution
$(I, \iota)$.  Then there is a canonical topology on each space $Ia$, with
the property that all the maps $If$ are continuous and closed and all the
maps $\iota_a$ are continuous.  Again, the topology can be defined in a
completely explicit way.

\begin{thm}	\label{thm:topo-terminal}
$(I, \iota)$ with this topology is the universal topological solution.
\end{thm}

Call a space \demph{self-similar} if it is homeomorphic to $Ia$ for some
self-similarity system $(\scat{A}, M)$ and some $a \in \scat{A}$, where
$(I, \iota)$ is the universal solution of $(\scat{A}, M)$.  There is a
`recognition theorem' giving a practical way to recognize universal
solutions, and this generates some examples of self-similar spaces:
\begin{itemize}
\item $[0,1]$, as in the Freyd example
\item $[0,1]^n$ for any $n\in\nat$; more generally, the product of two
  self-similar spaces is self-similar
\item the $n$-simplex $\Delta^n$ for any $n\in\nat$, by barycentric
  subdivision 
\item the Cantor set (isomorphic to two disjoint copies of itself)
\item Sierpi\'nski's gasket, and many other spaces defined by iterated
  function systems.
\end{itemize}

The proof of Theorem~\ref{thm:topo-terminal} involves showing that each of
the spaces $Ia$ is compact and metrizable (or equivalently, compact
Hausdorff with a countable basis of open sets).  So every self-similar
space is compact and metrizable.  The shock is that the converse holds:
\begin{thm}	\label{thm:ss-cm}
For topological spaces,
\[
\textrm{self-similar}
\iff
\textrm{compact metrizable}.
\] 
\end{thm}
This looks like madness, so let me explain.  

First, the result is non-trivial: the classical result that every nonempty
compact metrizable space is a continuous image of the Cantor set can be
derived as a corollary.  

Second, the word `self-similar' is problematic (even putting aside the
obvious objection: what could be more similar to a thing than itself?)
When we formalized the idea of a system of self-similarity equations, we
allowed ourselves to have infinitely many equations, even though each
individual equation could involve only finitely many spaces.  So there
might be infinite regress: for instance, $X_1$ could be described as a copy
of itself glued to a copy of $X_2$, $X_2$ as a copy of itself glued to a
copy of $X_3$, and so on.  Perhaps `recursively realizable' would be better
than `self-similar'.

Finally, this theorem does not exhaust the subject: it characterizes the
spaces admitting \emph{at least one} self-similarity structure, but a space
may be self-similar in several essentially different ways.

There's a restricted version of Theorem~\ref{thm:ss-cm}.  Call a space
\demph{discretely self-similar} if it is homeomorphic to one of the spaces
$Ia$ arising from a self-similarity system $(\scat{A}, M)$ in which the
category $\scat{A}$ is discrete (has no arrows except for identities).  The
Cantor set is an example: we can take $\scat{A}$ to be the one-object
discrete category.
\begin{thm}	
For topological spaces,
\[
\textrm{discretely self-similar}
\iff
\textrm{totally disconnected compact metrizable}.
\] 
\end{thm}
Totally disconnected compact Hausdorff spaces are the same as profinite
spaces, and the metrizable ones are those that can be written as the limit
of a \emph{countable} system of finite discrete spaces.  For instance, the
underlying space of the absolute Galois group
$\textrm{Gal}(\overline{\mathbb{Q}}/\mathbb{Q})$ is discretely
self-similar.

If you find the general notion of self-similarity too inclusive, you may
prefer to restrict to finite categories $\scat{A}$, which gives the notion
of \demph{finite self-similarity}.  This means that the system of equations
is finite.  A simple cardinality argument shows that almost all compact
metrizable spaces are not finitely self-similar.

I'll finish with two conjectures.  They both say that certain types of
compact metrizable space \emph{are} finitely self-similar.
\begin{conj}
Every finite simplicial complex is finitely self-similar.
\end{conj}
I strongly believe this to be true.  The standard simplices $\Delta^n$ are
finitely self-similar, and if we glue a finite number of them along faces
then the result should be finitely self-similar too.  For example, by
gluing two intervals together we find that the circle is finitely
self-similar.  Note that any manifold is as locally self-similar as could
be, in the sense of the introduction: every point is locally isomorphic to
every other point.

More tentatively,
\begin{conj}
The Julia set $J(f)$ of any complex rational function $f$ is finitely
self-similar. 
\end{conj}
This brings us full circle: it says that in the first example, we could
have taken any rational function $f$ and seen the same type of behaviour:
after a finite number of decompositions, no more new spaces $I_n$ appear.
Both $J(f)$ and its complement are invariant under $f$, so $f$ restricts to
an endomorphism of $J(f)$ and this endomorphism is, with finitely many
exceptions, a $\deg(f)$-to-one mapping.  This suggests that $f$ itself
should provide the global self-similarity structure of $J(f)$, and that if
$(\scat{A}, M)$ is the corresponding self-similarity system then the sizes
of $\scat{A}$ and $M$ should be bounded in terms of the degree of $f$.

\paragraph*{Acknowledgements}
I thank those who have given me the opportunity to speak on this:
Francis Borceux, Robin Chapman, Eugenia Cheng, Sjoerd Crans, Iain
Gordon, John Greenlees, Jesper Grodal, Peter May, Carlos Simpson,
and Bertrand To\"en.  I am very grateful to Jon Nimmo for
creating Figure~\ref{fig:Julia}(a).  I gratefully acknowledge a
Nuffield Foundation award NUF-NAL 04.

\end{document}